\documentclass[11pt]{article}

\usepackage{amssymb,amsfonts,amsmath,amsthm}
\usepackage{latexsym}
\usepackage{color}
\usepackage[hyperindex,colorlinks]{hyperref}

\newtheorem{thm}{Theorem}

\newtheorem{rem}{Remark}
\newtheorem{cor}{Corollary}

\newcommand{\setR}{{\mathord{\mathbb R}}}

\newcommand{\setH}{{\mathord{\mathbb H}}}

\newcommand{\setE}{{\mathord{\mathbb E}}}
\newcommand{\brf}[1]{\left\{{#1}\right\}}
\newcommand{\brr}[1]{\left({#1}\right)} 

\begin{document}

\title{On the composition of  finite rotations in $\setE^4$}

\author{Alex Goldvard\thanks{\small{Department of Mathematics, ORT 
Braude College, P.O. Box 78, 2161002 Karmiel, Israel, 
 \texttt{email:{goldvard@braude.ac.il} }}} \ and Lavi 
 Karp\thanks{\small{Department of Mathematics, ORT 
 Braude College, P.O. Box 78, 2161002 Karmiel, Israel, 
  \texttt{email:{karp@braude.ac.il} }}}}

\date{}

\maketitle

\begin{abstract}
We achieve compositions rules for the geometric parameters of the composed 
rotations, which is in a certain sense analogous to the well known Rodrigues 
formula.     
We also obtain  a necessary and sufficient  condition 
for a composition of two simple rotations in $\setE^4$ to be a simple 
rotation. 
\end{abstract}

\maketitle

\section{Introduction}

In this note we are studying  compositions of finite rotations in  the 
Euclidean space 
$\setE^4$, the four dimensional real vector space with the standard scalar 
product.  

There are several significant differences between rotations in 
$\setE^4$ and  in the  three dimensional space. The group of rotations 
in $\setE^3$ essentially  comprises  one type, that is, a rotation about an 
axis, while in $\setE^4$ 
there are three types of rotations: (i) \textit{A simple 
rotation}, 
it leaves  a plane (two dimensional subspace of $\setE^4$) point-wise fixed 
and  induces a  two dimensional rotation in the orthogonal plane; 
(ii) \textit{A Clifford translation 
(an isoclinic rotation)},  here each vector in $\setE^4$ turns through the same 
angle;  (iii) 
\textit{A double rotation}, here $\setE^4$ is decomposed into two 
orthogonal 
planes and points in the first plane rotate through an angle $\alpha$, 
while points in the second plane 
rotate through {an angle}, $\beta\neq\alpha$.

In this note we derive two properties of the compositions of rotations in 
$\setE^4$. The first one gives  analytical and geometrical 
characterizations of the  subgroup of  simple rotations. The second 
result  deals with  compositions formulas for double rotations, that is, 
we obtained formulas that enable the calculations of the orthogonal 
planes and angles  of the composed rotation  in 
terms of the corresponding  characteristics of each one of the rotations 
in $\setE^4$.

Our main tool is the quaternionic representation 
\begin{equation}
\label{MF} 
x\mapsto axb
\end{equation}
of rotations in $\setE^4$. 
Cayley obtained  formula  (\ref{MF}) in 1855  by means of laborious 
calculations \cite{Cayley}. 
 Many years later, in 1945, Coxeter~\cite{Coxeter}  gave a very elegant 
 proof of the above formula by employing   compositions of reflections 
 to represent rotations. He also derived the geometric characterization 
 of 
 the rotation~(\ref{MF}), that is, he computed the orthogonal planes and 
 angles via 
 the quaternions $a$ and $b$. But there are several cases in which 
 Coxter's formulas are not valid.  We shall fill the  gap
 in this note and complete   the  geometric characterization in 
 those cases.

The outline of the paper is as follows. In the next section we shall  give a 
brief account of Coxeter's approach and treat the above mentioned special cases. Section \ref{Main} 
 deals with the composition of simple rotations. In section \ref{GB} we 
 calculate planes and angles of a composed rotation and exhibit an example. We 
 end with section \ref{END} where we discuss the results.

\section{Representation of rotations and reflections by quaternions}

We first fix the notations. 
 A quaternion $x$ is expressed as $x=x_0+x_1i+x_2j+x_3k\equiv Sx+Vx$, where
 $Sx=x_0$ is called the real part, $Vx= x_1i+x_2j+x_3k$ is the vector part, and 
 the coefficients $x_0, x_1, x_2$  and $x_4$  
 are real numbers.  The  quaternions multiplication is defined by the table
 \begin{equation*}
 i^2=j^2=k^2=-1, \; ij=-ji=k, jk=-kj=i, ki=-ik=j.
 \end{equation*}
We denote by $\setH$ the set of all quaternions. 
Since a quaternion is completely 
 defined by its four real components, it determines a point (vector) in  the 
 four dimensional real space  $\setR^4$. Similarly,  $Vx=x_1i+x_2j+x_3k$ can be 
 identified with a point (vector) in three dimensional space $\setR^3$.  The 
 multiplication formula  
 \begin{equation}
 \label{msv}
xy=SxSy-Vx\cdot Vy+SxVy+SyVx+Vx\times Vy,
 \end{equation}
 is a consequence of the  distributive property of the quanternionic algebra. 
 Here 
 \\
$Vx\cdot Vy$ and $Vx\times Vy$ are  the scalar cross products in $\setE^3$ 
respectively.

Let $x$ be a quaternion, then $\bar x=Sx-Vx$ is the conjugate  quaternion, 
and    $Nx=\bar x x=x\bar x=x_0^2+x_1^2+x_2^2+x_3^2$ is the square of 
the norm. When $Nx=1$, then it is called a unit quaternion.
Any unit quaternion $a$ can be written in a polar form $a=\cos \alpha +p 
\sin 
\alpha$, where $p=p_1i+p_2j+p_3 k$ and  satisfies $p^2=-1$. The vector 
$p$ is 
called a pure unit  quaternion. A subspace which is spanned by two 
quaternions $a$ and $b$ is denoted by ${\rm Sp}\{a,b\}=\{s_1 a+s_2 b: s_1,s_2\in 
\setR\}$. 

The advantage of quaternionic representation  relies on two simple 
properties of quaternions. The function $f(x)=ax$ is linear and 
$N(ax)=NaNx$. Therefore, if $a$ is a unit 
quaternion,  then the multiplication $ax$ represents an 
orthogonal transformation in $\setE^4$.  
Moreover, for any unit quaternion $x$,  $(ax)\cdot 
x=S(ax\bar{x})=S(a)$, hence the cosine of 
the angle between $ax$ and $x$ is equal to $S(a)$ for any $x\in \setE^4$.
Thus any vector in $\setE^4$ rotates with angle $\alpha$, where 
$\cos\alpha=S(a)$.   Such type of rotations does not exist in 
$\setE^3$, and it is called a left Clifford translation. It is evident 
that the set of all left Clifford translations forms a subgroup of the 
rotations' group in $\setE^4$. Similarly, one may consider a right 
Clifford transformation $x\mapsto xb$, where $b$ is a unit quaternion.
Now, the composition of left and right Clifford translations results in 
the formula (\ref{MF}). As we already mentioned, Coxeter proved that 
this formula  comprises all rotations in $\setE^4$ by using reflections 
with respect to hyperplanes. For readers' convenient we recall his 
description of reflections and  rotations.

\subsection*{Reflections:}  Notice that any two quaternions $y$ and $x$ are 
orthogonal in $\setE^4$
if and only if $x\bar{y}+y\bar{x}=0$. Hence if $Ny=1$, then the transformation 
\begin{equation*}
\label{eq:reflection}
x\mapsto -y\bar{x}y
\end{equation*}
is a  reflection with respect to a hyperplane having $y$ as its normal (see 
\cite[\S 5]{Coxeter}). In what follows we denote this reflection by $R_y$.

\subsection*{Simple rotations:} A rotation in the plane  through the angle 
$\alpha$ and all the points in the orthogonal plane  remain fixed. Its 
 representation is given  by 
\begin{equation}
\label{eq:simple}
x\mapsto axb,\quad  Na=Nb=1,\quad Sa=Sb=\cos\frac{\alpha}{2}.
\end{equation}
It is shown in \cite{Coxeter} that
there are two unit quaternions $y$ and $z$ such that 
\begin{math}
	a=z\bar y, b=\bar yz
	\end{math},
	and the simple rotation (\ref{eq:simple}) is the composition of the 
	reflections 
	$ R_y(x)= -y\bar{x} y$ and $R_z(x) =-z\bar{x}z$, that is,   
\begin{equation}
\label{yz}
x\mapsto axb=(z\bar{y}) x(\bar{y} z).
\end{equation}
Obviously, $x\mapsto (z\bar{y}) x(\bar{y} z)=x$, whenever $x$ is orthogonal both 
to 
$y$ and $z$.  Hence,  the point--wise invariant plane of the simple rotation  
(\ref{yz})
is the intersection of the hyperplanes
\begin{equation*}
\Pi_1:=\{x\in\setH; \sum_{\mu=0}^3 y_\mu x_\mu =0\}, 
\Pi_2:=\{x\in\setH; \sum_{\mu=0}^3 z_\mu x_\mu =0\}.
\end{equation*}

\subsection*{Clifford translations:} 
Each vector in $\setE^4$ turns  through the same angle. 
Let $a$ be a unit quaternion, then $x \mapsto ax$ is a left 
Clifford translation and $x\mapsto xa$ is a right Clifford translation. 
The angle $\alpha$ of both types  satisfies the equation $\cos\alpha=Sa$.

Writing $a$ in a polar form $a=\cos \alpha+p\sin\alpha$, and consifer 
the left translation $x\mapsto ax$, then 
\begin{equation}
\label{Clifford}
x\mapsto x\cos\alpha +px\sin\alpha, \quad px\mapsto 
px\cos\alpha-x\sin\alpha.
\end{equation}
Hence the plane  ${\rm Sp}\{x,px\}$ is an invariant subspace of 
this transformation, and consequently  Clifford translations have 
infinite number of invariant two dimensional  subspaces. Similarly, for 
the 
right Clifford translation $x\mapsto xb$, $b=\cos\beta+q\sin\beta$, the 
planes ${\rm Sp}\{x,xq\}$ are invariant subspaces.  These subspaces are 
completely defined by the pure quaternions  $p$ and $q$ of the Clifford 
translations. In addition, we observe the similarity 
between~(\ref{Clifford})  
and the  formula for rotations in $\setE^2$.

\subsection*{Double rotations:} Two  orthogonal planes rotate simultaneously 
with 
two different angles. Let $p$ and $q$ be two pure quaternions such that 
$p^2=q^2=-1$, then a double rotation has the   representation
\begin{equation}
\label{eq:double}
x\mapsto (\cos\alpha+ p\sin\alpha)x(\cos\beta+ q\sin\beta)=:f(x).
\end{equation}
The rotations'  angles are  $\alpha\pm\beta$. In case that
$\cos\alpha=\cos\beta$, then the formula (\ref{eq:double}) represents a 
simple rotation.

Note that if $\cos\alpha, \cos\beta\neq 0$, then we can set 
\begin{math}
\tilde{p}=p\tan\alpha\end{math}, $\tilde{q}=q\tan\beta $
and rewrite formula (\ref{eq:double}) as
\begin{equation}
\label{gib-rep} 
f(x)=\cos\alpha(1+ \tilde{p})x(1+ \tilde{q})\cos\beta. 
\end{equation}
We will use this representation in Section~\ref{GB}. The quaternions 
$\tilde{p}$ and $\tilde{q}$ are called Gibbs vectors. We shall call 
(\ref{gib-rep}) Gibbs representation of the rotation. 

\subsection*{The orthogonal planes of the double and simple rotations:} 
Coxeter proved 
that the orthogonal 
planes  of the rotation (\ref{eq:double}) are 
\begin{equation*}
\Pi_1={\rm Sp}\brf{p-q,1+pq}\quad \text{and} \quad\Pi_2={\rm 
Sp}\brf{p+q,1-pq}.
\end{equation*}
Points in the plane $\Pi_1$  rotate through the angle 
$\alpha+\beta$  and in the plane $\Pi_2$ through the angle 
$\alpha-\beta$,  see  \cite[\S 9]{Coxeter}. It follows that if $\beta=\alpha$, 
then~(\ref{eq:double}) is a simple rotation with the points-wise invariant 
plane $\Pi_2$, while $\Pi_1$ is the fixed-plane when 
$\beta=-\alpha$.

If one of the  quaternions $p\mp q,1\pm q$ is zero, then  Coxeter's 
characterization of the invariant planes is not longer valid and we have 
to treat 
these cases separately. First of all note that the equalities $p=\pm q$ and 
$pq=\mp 1$ are equivalent and, in fact, there are  only two special cases: 
$q=\pm 
p$. It is readily follows from  formula~(\ref{eq:double}) that in these cases
\begin{equation*}
f(1)\in{\rm Sp}\brf{1,p} \quad \text{and} \quad f(p)\in{\rm Sp}\brf{1,p}.
\end{equation*}
Therefore $\Lambda_1={\rm Sp}\brf{1,p}$ is one of the rotating planes. 
Since the 
second one $\Lambda_2$, is the orthogonal complement of $\Lambda_1$, we 
conclude that
\[ \Lambda_2=\brf{x\in\setH; p\cdot x= Sx=0} \]
The pointwise invariant planes of  simple rotations when $q=\pm p$ can be found 
by the following considerations.
If $q=p$, then  by $f(1)\neq 1$,  and therefore $\Lambda_2$ is the fixed plane. 
While if $q=-p$, then $f(1)=1$, and hence $\Lambda_1$ is the fixed plane.

\section{Composition of  simple rotations in $\setE^4$}
\label{Main}
Let $f:x\mapsto axb$ and $g:x\mapsto cxd$ be two simple rotations and consider 
the 
composition $g\circ 
f: x\mapsto(ac)x(bd)$. We ask under which conditions $g\circ f$ will be again a 
simple 
rotation.

It follows from (\ref{eq:simple}) that the rotation $g\circ f$ is simple if and 
only if $S(ca)=S(bd)$. According to the multiplication formula (\ref{msv}), 
this condition is equivalent to 
\begin{equation*}
ScSa-Vc\cdot Va=SbSd-Vb\cdot Vc.
\end{equation*} 
Since both $f$ and $g$ are simple rotations, $S(a)=S(b)$ and $S(c)=S(d)$, hence 
the above condition results in  
\begin{equation}
\label{Vq}
Vc\cdot Va-Vb\cdot Vd=0.
\end{equation} 
We now recall that each rotation is decomposed into  two reflections by  
formula  (\ref{yz}). Hence, there  exists four unit quaternions 
$y,z,u,w$ such 
that 
\[a=z\bar y, b=\bar y z, c=w\bar u, d=\bar u w. \]
Hence
\begin{equation*}
\label{eq8}
\begin{split}
& Vc\cdot Va=V(w\bar u)\cdot V(z\bar y)  \\  = &
\left[-w_0Vu+u_0Vw+Vu\times Vw\right] \cdot 
\left[-z_0Vy+y_0Vz+Vy\times Vz\right]
\end{split}
\end{equation*}
and
\begin{equation*}
\label{eq9}
\begin{split}
& Vb\cdot Vd=V(\bar y z)\cdot V(\bar u w)  \\  = &
\left[y_0Vz-z_0Vy-Vy\times Vz\right] \cdot 
\left[u_0Vw-w_0Vu-Vu\times Vw\right].
\end{split}
\end{equation*}
Thus  (\ref{Vq}) becomes   
\begin{equation}
\begin{split}
Vc\cdot Va-Vb\cdot Vd &= -2\left\{y_0Vz\cdot\left(Vu\times Vw\right)-
z_0Vy\cdot\left(Vu\times Vw\right)\right. \\  &
\left. +u_0Vw\cdot\left(Vy\times Vz\right)-
w_0Vu\cdot\left(Vy\times Vz\right)\right\} \\ &=-2\det [y,z,u,w]=0,
\end{split}
\end{equation}
where $[y,z,u,w]$ denote  the matrix with the column vectors $y,z,u$ and  $w$. 
These vectors are the normals of the hyperplanes of the reflections $R_y, R_z, 
R_u$ and $R_w$ that decompose  the rotations. 
 Thus  we have obtained the following theorem. 

\begin{thm}
\label{thm1}
Let $f:x\mapsto axb$ and $g:x\mapsto cxd$ be two simple rotations such that $f$ 
is the 
composition of the reflections $R_y$ with $R_z$, and  $g$ is the composition of 
the reflections $R_u$ with $R_w$. Then the composition $g\circ f$  is a simple 
rotation if and only if  the four normals of the reflections are linearly 
dependent.
\end{thm}

The point-wise invariant planes of the rotations $f$ 
and $g$ are
\begin{equation*}
\Pi_f:=\{x\in\setH; \sum_{\mu=0}^3 y_\mu x_\mu=\sum_{\mu=0}^3 z_\mu 
x_\mu=0\}
\end{equation*}
and 
\begin{equation*}
\Pi_g:=\{x\in\setH; \sum_{\mu=0}^3 u_\mu x_\mu=\sum_{\mu=0}^3 
w_\mu x_\mu=0\}.
\end{equation*}

Obviously, if  $\Pi_f\cap\Pi_g\neq\{0\}$, then the composition 
$g\circ f$ is simple. However, the converse implication is not trivial. It 
follows from the above 
Theorem, that if the composition is simple, then the four normals are linearly 
dependent, which implies that   the intersection $\Pi_f\cap\Pi_g$
contains a  non--zero vector.  So a consequence  of the Theorem \ref{thm1} is: 

\begin{cor}
\label{cor}
Let $f:x\mapsto axb$ and $g:x\mapsto cxd$ be two simple rotations, and  $\Pi_f$ 
and 
$\Pi_g$ be 
the point--wise invariant planes of the rotations $f$ and $g$ respectively. 
Then 
the composition $h=g\circ f$ is a simple rotation if and only if $\Pi_f\cap 
\Pi_g \neq\{ 0\} $. 
\end{cor}

\begin{rem}
{ When the paper was written the authors were not aware that 
in 1890 Cole proved Corollary  \ref{cor}. 
 Cole essentially  worked out Cayley's parameters of the subgroup and 
 this method demands laborious computations. He also obtained formulas 
 for the orthogonal planes and the angle of the composed rotation, 
but they are complicated and lengthy, and inapplicable even  in some 
very simple examples.}
\end{rem}

\section{The angles and planes of the composed rotation}
\label{GB}

In this section we consider the following problem: Does there exist 
formulas that link the orthogonal planes and  angles  of two 
 rotations to those of composed rotation. 
 We start with double and simple rotations. So 
suppose that 
\begin{gather}
\label{eq:9}
f(x)= axb=\cos\alpha_1(1+\tilde p_1)x(1+\tilde q_1)\cos\beta_1\\
\label{eq:10}
g(x)= cxd=\cos\alpha_2(1+\tilde p_2)x(1+\tilde q_2)\cos\beta_2
\end{gather}
are two double rotations such that
\begin{equation}
\label{cond1}
\cos\alpha_1, \cos\beta_1\neq 0\quad \text{and } \ \cos\alpha_2, 
\cos\beta_2\neq 0.  
\end{equation}  
Recall that right hand side of (\ref{eq:9}) and (\ref{eq:10}) are Gibbs 
representation of the rotations $f$ and $g$.  
Consider the composition $f$ followed by $g$
\begin{equation*}
h(x):= (g\circ f)(x)=(ca)x(bd)=\cos\alpha(1+\tilde{p})x(1+\tilde{q})\cos\beta.
\end{equation*}
Then  in accordance to the multiplication formula 
 (\ref{msv}), 
 \begin{equation}
 \label{eq:11}
 \cos\alpha=S(ca)=\cos\alpha_1\cos\alpha_2(1-\tilde p_1\cdot \tilde p_2)
 \end{equation}
 and
 \begin{equation*}
 p\sin\alpha=V(ca)=\cos\alpha_1\cos\alpha_2(\tilde p_2+\tilde p_1 +\tilde p_2\times \tilde p_1).
 \end{equation*}  
Assuming
\begin{equation}\label{cond2}
\cos\alpha\neq 0,
\end{equation} 
we obtain that
\begin{equation}
\label{eq:12}
p\tan\alpha=\tilde{p}=\frac{\tilde p_2+\tilde p_1+\tilde p_2\times \tilde p_1}{1-\tilde 
p_2\cdot\tilde p_1}.
\end{equation}

In a similar manner,  
\begin{equation}
\label{eq:13}
\cos\beta=S(bd)=\cos\beta_1\cos\beta_2(1-\tilde q_1\cdot \tilde q_2),
\end{equation}
\begin{equation*}
q\sin\beta=V(bd)=\cos\beta_1\cos\beta_2(\tilde q_1+\tilde q_2+\tilde 
q_1\times \tilde q_2),
\end{equation*} 
and if
\begin{equation}\label{cond3}
\cos\beta\neq 0,
\end{equation} 
then  
\begin{equation}
\label{eq:14}
q\tan\beta=\tilde q=\frac{\tilde q_1+\tilde q_2+\tilde q_1\times \tilde q_2}{1-\tilde 
q_1\cdot\tilde q_2}.
\end{equation}
Thus in a view of those computations, we have 
obtained  the following theorem.

\begin{thm}
\label{PA}
If the conditions (\ref{cond1}), (\ref{cond2}) and (\ref{cond3}) hold, 
then the composition of the rotations (\ref{eq:9}) with (\ref{eq:10}), 
$h:=g\circ f$, is a rotation through the angles $\alpha\pm\beta$  about 
the planes 
spanned by quaternions $\tilde{p}\tan\alpha\pm\tilde{q}\tan\beta$ and 
$\tan\alpha \ \tan\beta\mp\tilde{p}\tilde{q}$, where 
$\alpha$, 
$\beta$, 
$\tilde{p}$ and $\tilde{q}$ are given by
(\ref{eq:11}) and (\ref{eq:13}), (\ref{eq:12}) and (\ref{eq:14}) respectively.  
\end{thm}

Note that the angels of the composed rotation are independent of the order of 
the  composition, that is, $g\circ f$ and $f\circ g$ result in the same angles. 
While the cross--product in (\ref{eq:12}) and (\ref{eq:14}) depend upon the 
order.

\begin{rem}
{ The expressions~(\ref{eq:12}) and (\ref{eq:14}) resemble the known 
Rodrigues formula  for the composition of rotations in 
$\setE^3$ (see e.g. \cite[ p.  58, 71]{Lounesto}). That is, if $f$ is a 
rotation through an angle $\alpha_1$ about an axis $c_1$, and $g$ is a rotation 
through an angle $\alpha_2$ about an axis $c_2$, then  their composition $g\circ 
f$ 
($f$ followed by $g$) is a rotation  through an angle $\alpha$ about an axis 
$c$, where   
\begin{equation}\label{Gibbs} 
\tilde c=\frac{\tilde c_2+\tilde c_1+\tilde c_2\times \tilde c_1}{1-\tilde 
c_2\cdot\tilde c_1}, 
\end{equation} 
 $\tilde{c}:=c\tan\frac{\alpha}{2},\tilde{c_1}:=c_1\tan\frac{\alpha_1}{2}$ and 
$ \tilde{c_2}:=c_2\tan\frac{\alpha_2}{2}$. 
The similarity of the formulas~(\ref{eq:12}),(\ref{eq:14}) 
and~(\ref{Gibbs}) is not accidental. A rotation in $\setE^3$ about an 
axis $c$ through an angle $\alpha$ has the quanternionic representation
$x\mapsto ax\bar{a}$, where $a=\cos\frac{\alpha}{2}+c\sin\frac{\alpha}{2}$. 
Therefore, the formula~(\ref{Gibbs}) is obtained in the same way as we 
do it for the composition of two double rotations in $\setE^4$. The 
vector $\tilde{c}$ is called  \textit{Gibbs vector} after J. W. Gibbs 
who invented it in 1881 in order to achieve formula~(\ref{Gibbs}). In 
fact this formula was obtained by Olinde  Rodrigues in 1840, 
see~\cite{Gray} for  historical survey.}
\end{rem}

\begin{rem}
{ Note that for a rotation in $\setE^3$, 
$\tilde{c}=c\tan\frac{\alpha}{2}$. Hence 
the discontinuity occurs 
when $\alpha=\pm \pi$, which means that the composed rotation is a reflection. 
For rotations in $\setE^4$ there are many possibilities for discontinuity. For 
example, in the case of a composition of simple rotations, where  
$\cos\alpha=0$ 
in (\ref{eq:11}), then it results in a reflection in one plane, and the 
identity 
in the orthogonal plan.
}\end{rem}

\begin{rem}
{ 
To the authors' knowledge Fedorov was the first one who discovered 
compositions rules of four dimensional rotations in 1958. 
Fedorov  interest was in applications of rotations to  relativity and 
quantum physics, hence he considered  Lorentz transformations in 
Minkowski space.  
He obtained  analogous formulas to (\ref{eq:11}) and (\ref{eq:13}), 
(\ref{eq:12}) and (\ref{eq:14}) in the Minkowski space and $\setE^4$. 
Fedorov showed that the group of 
 rotations in $\setE^4$ can be parametrized by two arbitrary three 
 dimensional vectors, and  by writing explicitly the matrices of 
 rotations in $\setE^4$ via those vectors he derived the compositions  
 formulas (\cite{Fedorov}, \S 24). So Fedorov  basically  used matrices 
 approach as opposed to the quaternions algebra. He did not 
 consider the geometrical characterization (orthogonal planes and 
 angles). Recently Fedorov's ideas on using Gibbs parametrization of a 
  rotation (it is also called a vector parametrization) became quite 
  widespread, see for example ~\cite{BMM},\cite{Muller} for further 
  discussions.
  
  }\end{rem}

We turn now to the composition of two Clifford translations of the same 
type.  So suppose $f:x\mapsto ax$ and $g:x\mapsto 
bx$ are two left Clifford translations, and consider the composition $f$ 
followed by 
$g$: $h=g\circ f: x\mapsto cx$, where $c=ba$.
Writing the quaternions by means of Gibbs vectors, that is, 
$a=\cos\alpha_1(1+\tilde{p}_1)$, $b=\cos\alpha_2(1+\tilde{p}_2)$ and 
$c=\cos\alpha(1+\tilde{p})$, then by the quanternionic algebra we get that 
\begin{equation*}
\tilde{p}=\frac{\tilde{p}_1+\tilde{p}_2+\tilde{p}_2\times 
\tilde{p}_1}{1-\tilde{p}_1\cdot\tilde{p}_2}
\end{equation*}
and 
\begin{equation*}
\cos\alpha=\cos\alpha_1\cos\alpha_2(1-\tilde{p}_1\cdot\tilde{p}_2).
\end{equation*}

Thus in view of formulas (\ref{eq:11}), (\ref{eq:12}), (\ref{eq:13}), 
(\ref{eq:14}), (\ref{Gibbs}) and the two last ones, we may say that
in a certain sense all type of rotations obey the same    composition 
rule.

\subsection{Example} 

The following example illustrates the consequences of Corollary \ref{cor} and 
Theorem  \ref{PA}.  
Let
 \begin{gather*}
 f:x\mapsto\frac{1+i}{\sqrt 2}x\frac{1+j}{\sqrt 2}\quad \text{and} \quad 
 g:x\mapsto\frac{1+j}{\sqrt 
 2}x\frac{1+k}{\sqrt 2}
 \end{gather*}
be  two simple rotations.
 Then $\Pi_f={\rm Sp}\brf{i-j,1+k}$ is the fixed--points plane of the 
 rotation $f$ 
 and
 $\Pi_g={\rm Sp}\brf{j-k,1+i}$ of $g$. We see that $\Pi_f\cap\Pi_g\neq 
 \{0\}$, hence 
 by Corollary \ref{cor} the composition $g\circ f$ should be a simple 
 rotation. Simple calculations show that  
 \begin{equation*}
 h=g\circ  f:x\mapsto\frac{1+i+j-k}{2}x\frac{1+i+j+k}{2}=:axb.
 \end{equation*}
Since $Sa=Sb$, the composition $g\circ f$ is simple.  
We shall now verify  formulas (\ref{eq:11}), (\ref{eq:12}) and (\ref{eq:14}).
 In order to it  we represent  these rotations by the  Gibbs vectors as in 
 (\ref{gib-rep}). So let  $\tilde{p}_f, \tilde{q}_f, \tilde{p}_g, \tilde{q}_g, 
 \tilde{p}_h, \tilde{q}_h$ the Gibbs vectors that
 corresponding those rotations, and let  $\alpha_f,\alpha_g,\alpha_h$ be
 angles of Gibbs vectors ($\tilde{p}= p\tan\alpha$, where $p$ is a pure 
 unit quaternion). Then
 
 \begin{gather*}
f:x\mapsto 
\left(\cos\frac{\pi}{4}\brr{1+i}\right)x 
\left(\brr{1+j}\cos\frac{\pi}{4}\right),\\
g:x\mapsto 
\left(\cos\frac{\pi}{4}\brr{1+j}\right)x \left(\brr{1+k}\cos\frac{\pi}{4}\right)
\end{gather*}
and
\begin{equation}
\label{polar}
h:x\mapsto \left(\cos\frac{\pi}{3}\brr{1+i+j-k}\right)x 
\left(\brr{1+i+j+k}\cos\frac{\pi}{3}\right).
\end{equation}
By (\ref{eq:11}),
\begin{gather*}
\cos{\alpha_h}=
\cos{\alpha_f}\cos{\alpha_g}\left(1-(i\cdot 
j)\right)=(\cos\frac{\pi}{4})^2=\frac{1}{2}.
\end{gather*}
Thus $\alpha_h=\frac{\pi}{3}$, which obviously coincides with (\ref{polar}).  
Now,  in accordance with formulas 
(\ref{eq:12}) and (\ref{eq:14}),
\begin{equation*}
\tilde p_h=\frac{j+i+j\times i}{1-j\cdot i}=i+j-k \quad\text{and} \ \ 
\tilde q_h=\frac{j+k+j\times k}{1-j\cdot k}=i+j+k.
\end{equation*}

\section{Conclusions}\label{END} 

There are various methods to represent rotations in $\setE^4$. In this article 
we are using the quaternionic algebra approach. It enables us to obtain the 
characterization  of a nontrivial subgroup of the simple rotations. Another 
advantage is the quaternionic  formula~(\ref{eq:double}) for double rotations 
that  describes in a very intrinsic way the geometric parameters. Relying upon 
that formula we had achieved the compositions rules by means of very 
simple calculations.

\end{document}